\documentclass[11pt,reqno,a4paper]{amsart}
\usepackage{amssymb}
\usepackage[cp1251]{inputenc}
\usepackage[russian,english]{babel}
\usepackage{graphicx}

\usepackage[colorlinks,pdfauthor=Yu.Brezhnev,
            bookmarksopen=false,bookmarks=false]{hyperref}

\allowdisplaybreaks[4]
\newcommand{\ri}{\mathrm{i}}
\newcommand{\re}{\mathrm{e}}
\newcommand{\ds}{\displaystyle}
\newcommand{\sss}{\scriptscriptstyle}
\newcommand{\ded}{{\boldsymbol\eta}}
\renewcommand{\,}{\mbox{\hspace{0.064em}}}
\newcommand{\ie}{\mbox{i.\hspace{0.091em}e.}}
\newcommand{\wpp}{{\textstyle\wp\mbox{\small$'$}}}
\newcommand{\Z}{\mbox{\large$\boldsymbol{\zeta}$}}
\newcommand{\WP}{\mbox{\large$\boldsymbol{\wp}$}}
\newcommand{\WPP}{\mbox{\large$\boldsymbol{\wp{\scriptstyle'}}$}}
\newcommand{\Btau}{\mbox{\Large$\boldsymbol\tau$}}
\newcommand{\dtheta}{{\ds\theta'{\!\!}_1}}
\newcommand{\deq}{\mbox{\rm \,\,\,\raisebox{0.034em}{:}$=$\,\,\,}}

\newcommand{\mfrac}[2] 
{\raisebox{0.047em}{\footnotesize$\displaystyle\frac{#1}{#2}$}}
\newcommand{\Mfrac}[2] 
{\raisebox{0.020em}{\small$\displaystyle\frac{#1}{#2}$}}

\newcommand{\bbig}[1]
{\raisebox{-0.1em}{\scalebox{1.1}[1.5]{\mbox{$#1$}}}}
\def\russian{\selectlanguage{russian}}

\begin{document}

\begin{abstract}
We represent and analyze the general solution of the sixth Painlev\'e transcendent in the Picard--Hitchin--Okamoto class in the Painlev\'e form as the logarithmic derivative of the ratio of  certain $\boldsymbol{\tau}$--functions. These functions are expressible explicitly in terms of the elliptic Legendre integrals and Jacobi $\theta$-functions, for which we write the general differentiation rules. We also establish a relation between the $\mathcal{P}_{6}$ equation and the uniformization of algebraic curves and present examples.
\end{abstract}

\author{Yu.\;V.\;Brezhnev}
\title[On \Btau-function solution to the $\mathcal{P}_6$-equation]
{A \Btau-function solution\\
to the sixth Painleve transcendent}
\thanks{Research supported by
the Federal Targeted Program under contract 02.740.11.0238.}

\hfill \noindent {\russian\small{\sc Теор.\ Мат.\ Физика} (2009), {\bf
161}(3), 346--366}

\hfill \noindent {\small{\sc Theor.\ Math.\ Phys.} (2009), {\bf
161}(3), 1616--1633}

\maketitle\thispagestyle{empty}\pagebreak

\tableofcontents

\section{Introduction}
\noindent The first appearance of the sixth Painlev\'e transcendent
\begin{equation}\label{P6}
\begin{aligned}
\mathcal{P}_{6}:\quad y_{\mathit{xx}}^{}={}&\frac12\!
\left(\frac1y+\frac{1}{y-1}+ \frac{1}{y-x}\right)\!
y_x^{}{}^{\!\!\!2}\,- \left(\frac1x+\frac{1}{x-1}+
\frac{1}{y-x}\right)\!y_x^{}\\[0.3em]
&{+}\,\frac{y(y-1)(y-x)}{x^2(x-1)^2}
\left(\alpha-\beta\,\frac{x}{y^2}+\gamma\,\frac{x-1}{(y-1)^2}-
\bbig(\delta-\mfrac12\bbig)\,\frac{x(x-1)}{(y-x)^2} \right)
\end{aligned}
\end{equation}
dates back to the 1905 Fuchs' paper \cite{fuchs}. Contemporary applications of the equation
are very well known
\cite{babich,tod,hitchin}  and the most nontrivial ones are in cosmology. Fuchs had already used the elliptic integral in his paper, and Painlev\'e a year later
wrote this equation in the remarkable form
\begin{equation}\label{P6wp} -\frac{\pi^2}{4}\,\frac{d^2
z}{d\tau^2}=\alpha\,\wpp(z|\tau)+
\beta\,\wpp(z-1|\tau)+\gamma\,\wpp(z-\tau|\tau)+
\delta\,\wpp(z-1-\tau|\tau)\,,
\end{equation}
again using elliptic functions and performing the transcendental variable change
$(y,x)\mapsto(z,\tau)$:
\begin{equation}\label{subs}
x=\frac{\vartheta_4^4(\tau)}{\vartheta_3^4(\tau)}\,,\qquad
y=\frac13+\frac13\frac{\vartheta_4^4(\tau)}{\vartheta_3^4(\tau)}
-\frac{4}{\pi^2} \frac{\wp(z|\tau)}{\vartheta_3^4(\tau)}
\end{equation}
(these formulae are equivalent to the original Painlev\'e ones; see relations \eqref{p} below). Here and hereafter, $\wp(z|\tau)\deq\wp(z|1,\tau)$ and other Weierstrass' functions
$(\sigma,\zeta,\wp,\wpp)$ \cite{a,bateman,we2} are constructed with the half-periods $(\omega,\omega')=(1,\tau)$, and the Jacobi theta constants $\vartheta(\tau)$ depend on the modulus $\tau$ belonging to the upper half-plane $\mathbb{H}^+$, \ie\
$\boldsymbol{\Im}(\tau)>0$ (see the appendix for definitions).

This Painlev\'e result, known as the $\wp$--form of the
$\mathcal{P}_{6}$-equation was
rediscovered and developed in the 1990s with other approaches, and the deep relation between   \eqref{P6} and
second-order linear differential equations of the Fuchsian class is well known. The transition from these
equations to Painlev\'e representation \eqref{P6wp} was recently well described informally \cite{bab}. In the same place additional references can also be found.

The interest in general Painlev\'e equation  \eqref{P6} increased dramatically after Tod demonstrated in 1994 \cite{tod}, that imposing the $\mathbb{SU}(2)$-invariance condition on the metric (the famous cosmological Bianchi-IX model) together with assuming the metric conformal invariance are consistent with the
Einstein--Weyl self-duality equations $R_{ab}=\frac12\Lambda g_{ab}^{}$ and
$W^+_{\mathit{abcd}}=0$. These equations are partial
differential equations but they can be reduced to an ordinary differential equation with respect to the
conformal factor, which can be expressed in terms of a solution of the
$\mathcal{P}_6$ equation for particular values of
the parameters
$(\alpha,\beta,\gamma,\delta)$ (sect.~2.1). Equation \eqref{P6} with arbitrary parameters was later obtained as
a self-similar reduction of the equations for the Ernst potential in general relativity \cite{schief}. This is not the
only example of the appearance of equation \eqref{P6} in applications: it suffices to recall the Yang--Mills equations (Mason--Woodhouse (1996)),
two-dimensional topological field theories (Dubrovin), etc.
In the general setting equation
\eqref{P6} was also derived in
\cite{nijhof}  as a self-similar reduction. See also references
in these works.

\subsection{The Painlev\'e substitution \eqref{subs}}
Equation \eqref{P6} and substitution \eqref{subs} reflect the properties of not
only the equation itself but also its solutions. If an equation is of the 2nd order, linear in  $y_{\mathit{xx}}^{}$, and
rational in $y_x^{}$ and if its solutions have only fixed branching points (the Painlev\'e property), then the number
of these points, as is known, is reducible to three. We can always choose them to be located at the points $x_j^{}=\{0,1,\infty\}$, and the equation then either takes Fuchs--Painlev\'e form \eqref{P6} or represents
some of its limit cases \cite{painleve}. In turn, the plane of the variable $x$  can be conformally mapped one-to-one to the
fundamental quadrangle of the modular group $\boldsymbol{\Gamma}(2)$ \cite{bateman} in the plane of the new variable
$\tau$ using the modular function $x=k'^2(\tau)$, \ie\ the first formula in \eqref{subs}. Function $k'^2(\tau)$ behaves exponentially with respect to the local parameter $\tau$ at the preimages of the points $x_j$. Because of this  the branching $y=y(x)$ of an arbitrary power-law or logarithmic character in the vicinity of $x_j$ transforms into a locally single-valued dependence on $\tau$ \cite[\textsection\,7]{br}. This partially explains\footnote{\label{ln}This is a partial explanation because a general rigorous statement on the branching of solutions of the $\mathcal{P}_{6}$ equation has
apparently not yet been introduced in the literature. For example, it is unknown whether it is (at most) a local branching
of the type of the rational function $x^\alpha\ln^n\! x$. See the important paper \cite{guz}, which was devoted to precisely this question and
where Fuchsian equation \eqref{leg}---an equivalent to equation \eqref{P6wp}---was the main object of study.} the origin of substitution \eqref{subs}. Fuchs obtained this substitution  \cite{fuchs} using
the Legendre equation and Painlev\'e used the modular function \cite{painleve}.

\footnotesize
Painlev\'e himself wrote the equation not in form \eqref{P6wp} but in terms of the modular function $x=\varphi(\xi)$ with the three singular points $x_j=\{0,1,\infty\}$, appearing when inverting the elliptic integral

\begin{equation}\label{p}
z=\frac{1}{\omega_1^{}}\int_\infty^{y}
\!\!\frac{dy}{\sqrt{y(y-1)(y-x)}}\,,\qquad
\frac{\omega_2^{}}{\omega_1^{}}=\xi\,,\qquad y=f(z,\xi)\,.
\end{equation}
In his original notation, the equation $\mathcal{P}_{6}$  in the $\wp$-form is as follows
\cite[p.\;1117]{painleve}:
\begin{figure}[htbp]
\centering
\includegraphics[width=13 cm,angle=0]{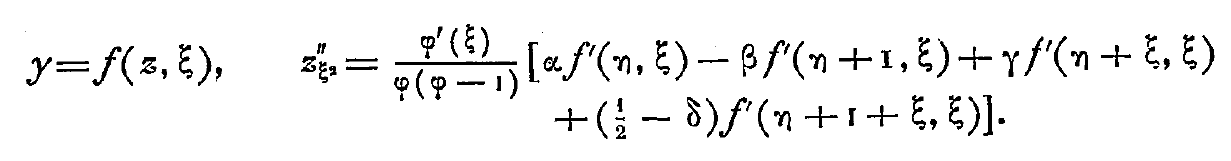}
\end{figure}

\noindent
where the misprint  $\eta\mapsto z$ should be corrected.
Equation \eqref{P6},
has the form of an `inhomogeneous' Legendre differential equation
\cite{fuchs}
\begin{equation}\label{leg}
\frac{d^2u}{dt^2}+\frac{2\,t-1}{t(t-1)}\,\frac{du}{dt}+\frac{u}{4\,t(t-1)}=
\frac{\sqrt{\lambda(\lambda-1)(\lambda-t)}}{2\,t^2(t-1)^2} \left[
k_\infty-k_0^{}\frac{t}{\lambda^2}+k_1^{}\frac{t-1}{(\lambda-1)^2}-
k_t^{}\frac{t(t-1)}{(\lambda-t)^2} \right],
\end{equation}
where
$u={\mbox{\raisebox{-0.2em}{${\scalebox{1.1}[1.5]{$\int$}}_{{\!\!\textstyle{}_0}}^\lambda$}}}
\!\frac{d\lambda}{\sqrt{\lambda(\lambda-1)(\lambda-t)}{}}$
and $\lambda=\lambda(t)$ is a solution of the equation  $\mathcal{P}_{6}$.
\normalsize

\subsection{The motivation for the work and the paper content}
Other singular points are movable but are poles,
and the global behavior of solutions is therefore totally governed by the logarithmic derivatives of `entire' functions
\cite{painleve}, \cite[pp.\;77--180]{conte}. Painlev\'e stressed that construction of such transcendental functions completes the procedure of integrating an equation (the concept of the `int\'egration parfaite' in the Painlev\'e terminology). The
pole character is known for all Painlev\'e equations \cite{gromak,conte} and all their solutions except those for the equation
$\mathcal{P}_{1}$ have the structure \cite[p.\;123]{painleve}
\begin{equation}\label{entire}
y\sim \frac{d}{dx}\mathrm{Ln}\, \frac{\Btau_{\!\!2}^{}}
{\Btau_{\!\!1}^{}}
\end{equation}
(in the generic case) with the functions $\Btau_{\!\!1,2}^{}(x)$ having no movable singularities. See the works \cite[p.\;371]{okamoto2} and \cite[p.\;165]{conte} for the complete list of such formulae.

If we disregard the known automorphisms of solutions of the equation
$\mathcal{P}_{6}$ in the space of the parameters
$(\alpha,\beta,\gamma,\delta)$ (see, e.g., the Gromak paper in \cite{conte}, \cite{korotkin}, \cite[\textsection\,42]{gromak}, \cite{iwasaki} and the references therein), only two
cases are known where the general integral for equation \eqref{P6}
can be written. These are the solutions of Picard
\cite[pp.\;298--300]{picard} and Hitchin (sect.\;2.1),
but Painlev\'e form \eqref{entire} is known for neither of these solutions.

On the other hand, representing solutions in terms of the corresponding analogues of $\Btau$-functions occurs
in almost all investigations of nonlinear equations somehow related to the integrability property: in the
Hirota method, in the theory of isomonodromic deformations \cite{korotkin}, in soliton theory, and in its $\Theta$-functional
generalizations. In this respect, it is interesting to show that the Painlev\'e equations also admit a nontrivial
situation in which the functions $\Btau_{\!\!1,2}^{}(x)$
can be written explicitly when we have the general integral of
motion. Here, we fill this gap (sects.\;2.2, 4) and also describe the distribution of solution
singularities (sect.\;3). Technically, we need special differential properties of Jacobi's theta-functions and their
theta-constants. These properties are also new, and we present them in the appendix. As follows from
the solution (this was properly noted in \cite{korotkin}), the character of the Painlev\'e
$\Btau$-functions differs drastically from that of the
$\Theta$-functional solutions of the solitonic equations. Nevertheless, a variant of the straight-line
section of a manifold is present (sect.\;4). In section~5 we demonstrate the connection of the equation $\mathcal{P}_{6}$ with the uniformization of nontrivial algebraic curves. To what was said in sect~1.1, we add
that Painlev\'e himself explicitly wrote \cite[p.\;81]{conte} about representing solutions in terms of single-valued functions in the context of the Painlev\'e property. Picard mentioned the same thing many times in \cite[pp.\;93, 188]{PO},
\cite{picard}. The relation to uniformization mentioned below and the presentation of numerous examples
(sect.~5.2) is therefore not accidental.

\section{The Painlev\'e form}

\subsection{The Hitchin solution} The Hitchin solution corresponds to the parameter values
\begin{equation}\label{ab}
\alpha=\beta=\gamma=\delta=\frac18\,,
\end{equation}
which were found in \cite{hitchin} under the condition that the Einstein equation metric is conformal with $\mathbb{SU}(2)$-invariant anti-self-dual Weyl tensor \cite{tod}. Equation \eqref{P6wp} in this case can be easily reduced to the equivalent form
written in theta functions
\begin{equation}\label{th}
\frac{d^2 z}{d\tau^2}=4\pi\,\ded^9(\tau)\,
\frac{\theta_1(2z|\tau)}{\theta_1^4(z|\tau)}\,,
\end{equation}
where $\ded$ is the Dedekind function
$$\ded(\tau)=\re^{\frac{\pi\ri}{12}\,\tau}_{\mathstrut}\,
{\prod\limits_{\sss k=1}^{\sss \infty}}^
{\ds\mathstrut}_{\ds\mathstrut}
\big(1-\re^{2k\pi\ri\,\tau}\big)\,.
$$
The general solution to equation \eqref{th} in parametric form \eqref{subs} is \cite{hitchin}
\begin{equation}\label{hitchin1}
\wp(z|\tau)=\wp(A\tau+B|\tau)+\frac12\, \frac{\wpp(A\tau+B|\tau)}
{\zeta(A\tau+B|\tau)-A\,\eta'(\tau)-B\,\eta(\tau)}\,,
\end{equation}
where $A,B$ are the integration constants. This solution in a different implicit form was mentioned by
Okamoto \cite[p.\;366]{okamoto}.

Solution \eqref{hitchin1} is not only the most nontrivial of all those currently known but is also rather instructive for $\mathcal{P}$-equations in general because it provides the general integral. It has been considered in detail in the literature
and is certainly more than `just a solution'. The Picard solution is degenerate in this respect because it is a
perfect square (see sect.~4.1). This explains why we choose the Hitchin solution for further analysis. Other
particular cases are known that contain variations of the Painlev\'e
$\partial_x\mbox{Ln}$-form but do not contain \Btau-functions. Allowing more freedom in choosing the parameters $(\alpha,\beta,\gamma,\delta)$, these solutions nevertheless contain only one integration constant and are related to the linear hypergeometric
${}_2F_1$-equation \cite[p.\;733]{conte},
\cite[\textsection\,44]{gromak}, \cite[p.\;374
\textsection\,5.4]{okamoto}, \cite[p.\;145]{iwasaki}.

\subsection{The Painlev\'e form}
One of the theta-function versions of solution \eqref{hitchin1} was proposed in
\cite{hitchin} and another form in \cite{korotkin}. But the complexity of these solution forms was mentioned in \cite[p.\;75]{hitchin} itself, in  \cite[p.\;901]{korotkin}, and also in work \cite{babich}, consisting entirely of calculations. In those papers, solutions were also presented in the
parametric forms and contained the set of functions $\theta,\theta',\theta'',\theta''',\vartheta,\vartheta',
\vartheta'''$. However one can obtain the
solution in a compact form if we transform Weierstrass functions  in \eqref{hitchin1} into the theta functions
(cf. \cite[p.\;74]{hitchin}).

{\bf Proposition 1.} {\em The general solution to equation  \eqref{P6} with parameters \eqref{ab} is}
\begin{equation}\label{simp}
y=\frac{\sqrt{x}}{\theta_1^2}
\left\{\frac{\pi\,\vartheta_2^2\!\cdot\!\theta_2\,\theta_3\,\theta_4}
{\dtheta+\pi\, \ri\,A\,\theta_1}-\theta_2^2\right\},
\end{equation}
{\em where $\theta=\theta\big(\frac12A\tau+\frac12B|\tau\big)$ and
$\vartheta=\vartheta(\tau)$.}

The possibility of representing solution \eqref{simp} in explicit form \eqref{entire}
is related to nontrivial properties of
the theta functions, for instance, to the identity that follows from Lemma\;3 in the appendix
$$
\frac{\pi}{2\,\ri}\,\frac{d}{d\tau}
\mathrm{Ln}\big\{\zeta(z|\tau)-z\,\eta(\tau)\big\}=
\wp(z|\tau)+\frac12\,
\frac{\wpp(z|\tau)}{\zeta(z|\tau)-z\,\eta(\tau)}+\eta(\tau)\,,
$$
Comparing this property with \eqref{hitchin1}, we obtain the total logarithmic derivative,
\begin{equation}\label{wpz}
\wp(z|\tau)=\frac{\pi}{2\,\ri}\,\frac{d}{d\tau}\mathrm{Ln}
\frac{\zeta(A\tau+B|\tau)-A\,\eta'(\tau)-B\,\eta(\tau)}
{\ded^2(\tau)}\,.
\end{equation}
Replacing $(A,B)\mapsto 2(A,B)$ and transforming the right-hand side of
\eqref{wpz} into $\theta$-functions (lemma~\;1 sect.~7) we can write the solution in the form
\begin{equation}\label{Tau}
y=\frac{2\,\ri}{\pi}\frac{1}{\vartheta_3^4(\tau)}\, \frac{d}{d\tau}
\mathrm{Ln}\frac{\dtheta(A\tau+B|\tau)+2\,\pi\,\ri\,A\,
\theta_1\!(A\tau+B|\tau)}
{\vartheta_2^2(\tau)\,\theta_1\!(A\tau+B|\tau)}\,, \qquad
x=\frac{\vartheta_4^4(\tau)}{\vartheta_3^4(\tau)}\,.
\end{equation}
To rewrite the solution in the initial variables $(x,y)$, we must use the formulae inverting substitution \eqref{subs} and its differential. They can be written in terms of the complete elliptic integrals
$K,\,K'$ \cite{a,WW}. It is
easy to see that this transition together with the differentiation is determined by the formulae
\begin{equation}\label{LnK}
\frac{d}{d\tau}=\pi\ri\,x\,(x-1)\,\vartheta_3^4(\tau)\,
\frac{d}{dx}\,,\qquad
\vartheta_2^2(\tau)=\frac{2}{\pi}\,\sqrt{1-x}\,K'(\sqrt{x})\,,
\end{equation}
and inversion of the first formula in \eqref{subs} is given by the expression (series)
\begin{eqnarray}
\ds
\re^{\pi\ri\tau}&\!\!\!\!=&\!\!\!\exp\Big\{{-}\pi\mfrac{K(\sqrt{x})}
{{\ds K'}(\sqrt{x})}\Big\}=\label{weier}
\\\nonumber\\\ds
\phantom{e^{\pi\ri\,\tau}} &\!\!\!\!=&\!\!\!
\Big(\mfrac{1-x}{16}\Big)+8 \Big(\mfrac{1-x}{16}\Big)^{\!2}
+84\Big(\mfrac{1-x}{16}\Big)^{\!3} +
992\Big(\mfrac{1-x}{16}\Big)^{\!4}
+12514\Big(\mfrac{1-x}{16}\Big)^{\!5}+\cdots\,. \nonumber
\end{eqnarray}
Simplifying the obtained expressions and replacing $\ri A\mapsto  A$, we obtain the sought solution form.

{\bf Theorem.} {\em The general integral of equation \eqref{P6} with parameters
\eqref{ab} has the form}
\begin{equation}\label{hitchin}
y=2\,x\,(1-x)\,\frac{d}{dx}\mathrm{Ln}
\mfrac{\mbox{\normalsize$\dtheta$}\!\bbig(
A\frac{K}{K'}+B\bbig|\frac{\ri\, K}{K'}\bbig)+
\mbox{\normalsize$2\pi A\!\cdot\!\theta_1$}\!\bbig(
A\frac{K}{K'}+B\bbig|\frac{\ri\, K}{K'}\bbig) }
{\mbox{\normalsize$\sqrt{1-x}\,K'\cdot\theta_1$} \!\bbig(
A\frac{K}{K'}+B\bbig|\frac{\ri\, K}{K'}\bbig)}\,.
\end{equation}
{\em Here and occasionally hereafter, we use the shorthand notation $K\deq K(\sqrt{x})$ and  $K'\deq
K'(\sqrt{x})$.}

The quantities $K,\,K'$ can be written in terms of the hypergeometric functions
\cite{a,bateman}
$$
K(\sqrt{x})=\frac{\pi}{2}\cdot
{}_2F_1\Big(\mfrac12,\mfrac12;1\bbig|x \Big)\,,\qquad
K'(\sqrt{x})=\frac{\pi}{2}\cdot
{}_2F_1\Big(\mfrac12,\mfrac12;1\bbig|1-x \Big)\,,
$$
but we can use these series simultaneously only in the intersection of the convergence domains. In applied
problems, we therefore avoid using numerous transformations for analytic continuations of the ${}_2F_1$-series by regarding the functions $K,\,K'$ as the complete elliptic integrals
$$
K(\sqrt{x})=\int\limits_0^{\,\,1}\!\!\frac{d\lambda}
{\sqrt{(1-\lambda^2)(1-x\,\lambda^2)}}\,,\qquad
K'(\sqrt{x})=\int\limits_x^{\,\,1}\!\!\frac{d\lambda}
{\sqrt{(\lambda^2-x)(1-\lambda^2)}}
$$
(which are everywhere finite functions) or as the Legendre functions
$P_{\!-\frac12},\,Q_{\!-\frac12} $\cite{WW}.

Verifying this by substituting \eqref{hitchin} directly in \eqref{P6} is a rather nontrivial exercise in which we must use Lemma\;1 and its
simple corollaries. We have to use also the known differentiation rules for complete elliptic integrals \cite{a}. Analogously setting $E\deq E(\sqrt{x})$ and
$E'\deq E'(\sqrt{x})$, we can write these rules as
\begin{equation}\label{EK}
\begin{array}{rlrl}
\ds 2\,\frac{d}{dx}K &\!\!\!\ds=\frac{E}{x\,(1-x)}-
\frac{K}{x}\,,&\ds 2\,\frac{d}{dx}K'&\!\!\!\ds=\frac{E'}{x\,(x-1)}-
\frac{K'}{x-1}\,,\\\\
\ds 2\,\frac{d}{dx}E&\!\!\!\ds=\frac{E}{x}- \frac{K}{x}\,, &\ds
\qquad 2\,\frac{d}{dx}E'&\!\!\!\ds=\frac{E'}{x-1}- \frac{K'}{x-1}
\end{array}
\end{equation}
and they imply the Legendre relation
$$
\eta'=\tau\,\eta-\frac{\pi}{2}\ri\quad \Leftrightarrow\quad
E\,K'+E'K-K\,K'=\frac{\pi}{2}\,.
$$

Speaking more rigorously, in  \eqref{hitchin}, we must cancel the common `non-entire' factor $\exp\big(\frac14\pi\ri\tau\big)$
in the series
for  $\theta_1$ and $\dtheta$. But this factor and also the nonmeromorphic singularity of the solution
determined by the factor
$K'(\sqrt{x})$ are `fixed'. Relations \eqref{hitchin} and \eqref{EK} demonstrate that the solution contains
an additive term, which has a fixed critical singularity of the logarithmic type.

{\bf Corollary 1.} {\em The function $y(x)$ and the solution $z(\tau)$ of equation \eqref{th} have the forms}
\begin{equation}\label{52}
y=\frac{E'}{\ds K'}+2\,x\,(1-x)\,\frac{d}{dx}\mbox{Ln}\mspace{-2mu}
\left\{\frac{\dtheta}{\theta_1}\bbig( \textstyle
A\frac{K}{K'}+B\bbig|\frac{\ri\, K}{K'}\bbig)+2\pi A\right\}\,,
\end{equation}
\begin{equation}\label{53}
z=\wp^{\mbox{-}1}\!\!\left(\frac{\pi}{2\ri}\,\frac{d}{d\tau}\mbox{Ln}
\!\left\{\Mfrac{\dtheta(A\tau+B|\tau)}{\theta_1(A\tau+B|\tau)}
+2\pi\ri A\right\}-\eta(\tau)\Big|\tau\right),
\end{equation}
{\em where the elliptic integral $\wp^{\mbox{-}1}$ here and hereafter is defined as}
$$
\wp^{\mbox{-}1}(u|\tau)\deq
\int\limits_\infty^{\,\,u}\!\!\frac{d\lambda}
{\sqrt{4\lambda^3-g_2^{}(\tau)\lambda-g_3^{}(\tau)}}\,.
$$

To follow the meromorphic $x$-component (the movable poles) more closely, we substitute expansions of
type \eqref{weier} in formulae \eqref{hitchin}, \eqref{52} and in the theta-function series in the appendix. Weierstrass derived
analogous series in relation to the solution of the modular inversion problem for the function $k^2(\tau)$ \cite[pp.\;53--4, 56, 58]{we2}, \cite[p.\;367]{WW}.
Explicit expansions in the vicinities of all other points can be written using the formulae in the appendix.
Formulae \eqref{Tau},
\eqref{hitchin}, \eqref{52} themselves provide an analytic answer for the power-law expansions in \cite{hitchin} on pages 89--92 and 108--9.

\section{The pole distribution}
\noindent One of the two series of poles $x_{\mathit{mn}}$ of solution
\eqref{hitchin} admits an explicit parameterization,
\begin{equation}\label{poles}
\theta_1 \!\!\left(\!
A\mfrac{K(\sqrt{x})}{K'(\sqrt{x})}+B\Big|\mfrac{\ri\,
K(\sqrt{x})}{K'(\sqrt{x})}\right)=0\quad\Rightarrow\quad
x_{\mathit{mn}}=
\frac{\vartheta_4^4}{\vartheta_3^4}\Big(\mfrac{m-B}{n+A}\Big)\,,
\quad n,m\in\mathbb{Z}\,.
\end{equation}
We must supplement it with the natural restriction on  $(n,m)$ of the form
$\boldsymbol\Im\bbig(\frac{m-B}{n+A} \bbig)>0$, which imposes no restrictions on $A$ and $B$  other than these numbers cannot be real
simultaneously. Although this pole distribution is described by the simple formula \eqref{poles}, its actual form is
rather involved (see  {\sc Fig.\;1}).

\begin{figure}[htbp]
\centering
\includegraphics[width=10 cm,angle=0]{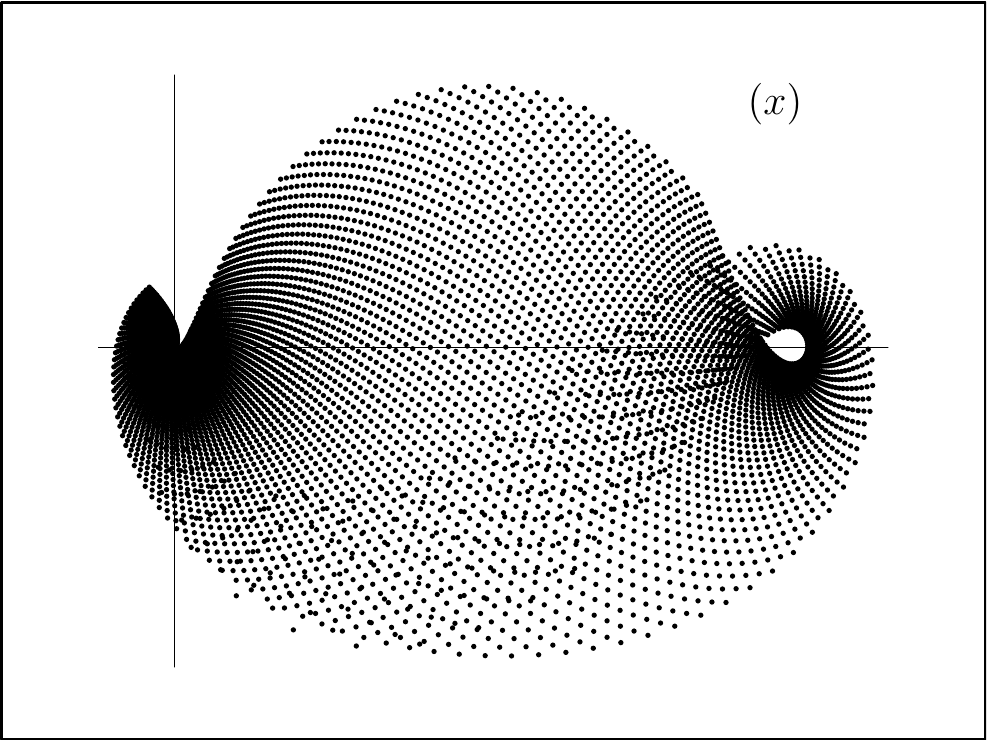}
\caption{The poles $x_{\mathit{mn}}$ of form \eqref{poles} for solutions
\eqref{hitchin}, \eqref{pic} under $A=125.45-103.29\,\ri$,
$B=36.710-69.980\,\ri$ and $(n,m)=-30\ldots 70$.}
\end{figure}
The general and characteristic behavior of the `orthogonal' directrices of the parameterizations  $n=\mbox{const}$ and $m=\mbox{const}$, \ie\ the pattern in
{\sc Fig.\;1} as a whole is shown in {\sc Fig.\;2}.
\begin{figure}[htbp]
\centering
\includegraphics[width=10 cm,angle=0]{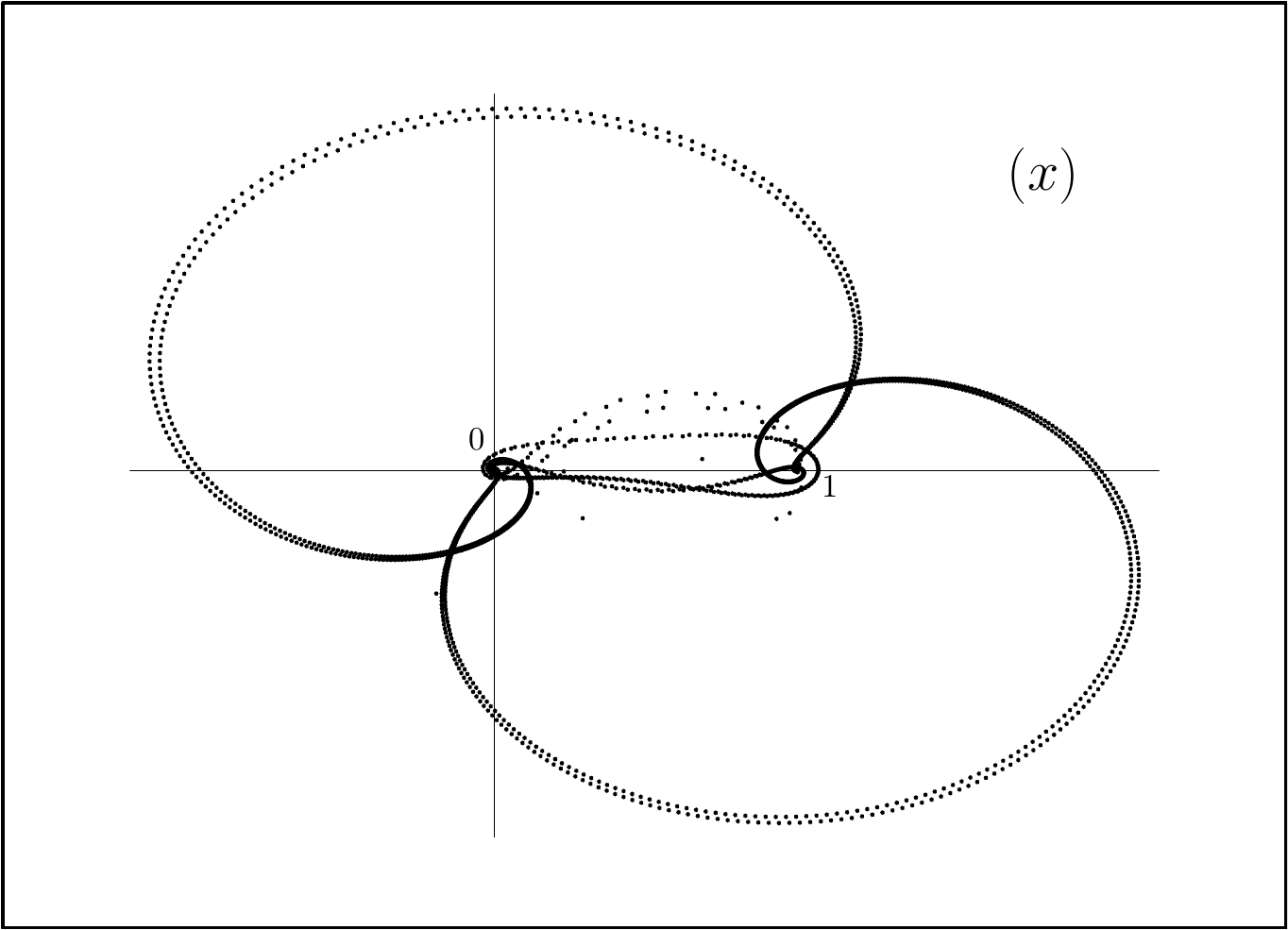}
\caption{Two pairs of directrices $n=(12,13)$ and $m=(129,130)$
of pole distributions \eqref{poles} for $A=-195.45+103.29\,\ri$,
$B=-6.710+79.98\,\ri$. All the seemingly `casual' points are also `good' poles belonging to these directries.}
\end{figure}
The pattern in this figure can be
interpreted as an analogue, although distant, of the two straight-line orbits of poles of any double periodic
function. The pole distributions are deformed when the initial data $(A,B)$ are changed, which generates
various patterns, but we observe that all the poles concentrate near the fixed singularities $x_j=\{0,1,\infty \}$,
as should be the case. We thus obtain the character of this deformation. If we take
one pole $p=x_{\mathit{mn}}$ under fixed $(m,n)$, then we find that it moves without changes along the straight-line
characteristic $(n+A)\,t-(m-B)=0$ in the space of the initial data $(A,B)$. The dynamics $p=p(t)$ for any pole of series \eqref{poles} is given by the same function
$p(t)=\frac{\vartheta_4^4}{\vartheta_3^4}(t)$ but with times
$t_k^{}$ that differ
by a fractional linear transformation for each pole $p_k^{}=p(t_k^{})$. The function $p(t)$ and the structure of the
fundamental domain of its automorphism group are well studied. Real or purely imaginary poles can also
be easily described.

We further note that the Okamoto transformation
\cite[\textsection\,47]{gromak}, \cite{korotkin}
\begin{equation}\label{ok}
\begin{array}{rl}
\ds\mbox{\small\sc Picard}_y\mapsto \mbox{\small\sc Hitchin}_{\tilde
y}:\quad&\ds \tilde
y=y+\frac{y\,(y-1)\,(y-x)}{x\,(x-1)\,y_x^{}-y^2+y_{{}_{\ds\mathstrut}}}\,,\\
\ds\mbox{\small\sc Hitchin}_y\mapsto \mbox{\small\sc Picard}_{\tilde
y}:\quad&\ds \tilde y=y-\frac{y\,(y-1)\,(y-x)}{
x\,(x-1)\,y_x^{}+\frac12\,y^2-x\,y+\frac12\,x}\,,
\end{array}
\end{equation}
which is known to transform the Picard and Hitchin solutions one into the other, leaves poles \eqref{poles} invariant
but creates or, under the inverse transformation, annihilates the second series of poles determined by
roots of the equation
$$
\zeta(A\tau+B|\tau)=A\,\zeta(\tau|\tau)+B\,\zeta(1|\tau)\,.
$$
Solving this equation is equivalent to finding the $A$-points of the transcendental function
\begin{equation}\label{mero}
f(x;A,B)=
\frac{1}{2\pi}\mfrac{\mbox{\normalsize$\dtheta$}\!\left(\!
A\frac{K(\sqrt{x})}{K'(\sqrt{x})}-B\bbig|\frac{\ri\,
K(\sqrt{x})}{K'(\sqrt{x})}\right)}{\mbox{\normalsize$\theta_1$}
\!\left(\! A\frac{K(\sqrt{x})}{K'(\sqrt{x})}-B\bbig|\frac{\ri\,
K(\sqrt{x})}{K'(\sqrt{x})}\right)}\,,
\end{equation}
\ie\ to solving the equation $f(x;A,B)=A$. The pole series determined by relation \eqref{mero} has a much more
involved description and evolution than for pole series \eqref{poles}.

\section{The \protect\Btau-functions}
\subsection{The Picard solution}

Pole series \eqref{poles} also exhausts all the poles of solutions to
equation \eqref{P6} in the Picard case $\alpha=\beta=\gamma=\delta=0$. It then follows from \eqref{P6wp} that $\ddot z=0\;\Rightarrow\;z=a\tau+b$.
Writing relations \eqref{subs} in terms
of the theta function, we obtain another form of substitution \eqref{subs} and the equivalent form of the Picard
solution in \cite{picard}:
\begin{equation}\label{pic}
y=-\frac{\vartheta_4^2(\tau)}{\vartheta_3^2(\tau)}\,
\frac{\theta_2^2\big(\frac z2|\tau \big)}{\theta_1^2\big(\frac
z2|\tau \big)}\qquad\Rightarrow\qquad
y_{{}_{\mbox{\tiny\sc Pic}}}=-\sqrt{x}\;
\mfrac{\mbox{\normalsize$\theta_2^2$}\!\left(\!
A\frac{K(\sqrt{x})}{K'(\sqrt{x})}+B\bbig|\frac{\ri\,
K(\sqrt{x})}{K'(\sqrt{x})}\right)}{\mbox{\normalsize$\theta_1^2$}
\!\left(\! A\frac{K(\sqrt{x})}{K'(\sqrt{x})}+B\bbig|\frac{\ri\,
K(\sqrt{x})}{K'(\sqrt{x})}\right)}\,.
\end{equation}

\footnotesize A comment is in order. Strictly speaking, the original Picard solution $u_{{}_{\mbox{\tiny\sc Pic}}}$ \cite[p.\;299]{picard} and the Fuchs--Painlev\'e
representation  \eqref{P6} are not identical but are related by an inverse transformation. Picard constructed
an example of the equation using the Legendre form $\Delta x=\sqrt{(1-x^2)(1-k^2x^2)}$, while Painlev\'e used
Riemann form \eqref{p}. The transition between solutions of these equations is therefore given by the formula
\begin{equation}\label{su}
u_{{}_{\mbox{\tiny\sc
Pic}}}\!\big({\textstyle\frac1x}\big)=2\,\frac{\sqrt{e-\wp}}{\pi\,\vartheta_3^2}\,=
\sqrt{y_{{}_{\mbox{\tiny\sc Pic}}}\!(x)}\,.
\end{equation}
In other words, we have the property that the case of parameters
$\alpha=\beta=\gamma=\delta=0$
is a case (not unique)
of Painlev\'e equation \eqref{P6} such that the square root of any of its solution satisfies, in turn, an equation that
also has the Painlev\'e property, namely, the Picard equation
\cite[p.\;299]{picard}
\begin{equation}\label{u}
\frac{d^2u}{dx^2}-\Big(\frac{du}{dx}
\Big)^{\!2}\,\frac{u\,(2xu^2-1-x)}{(1-u^2)(1-xu^2)}+
\frac{du}{dx}\!\left[
\frac{u^2-1}{(1-x)(1-xu^2)}+\frac1x\right]-
\frac{u\,(1-u^2)}{x\,(1-x)(1-xu^2)}=0\,,
\end{equation}
where we must amend the last term with the omitted factor
$\frac14$. This equation appears at the end of the
almost two-hundred-page-long Picard treatise  \cite[p.\;298]{picard} as \'equation diff\'erentielle
curieuse satisfied by the function
$u=\mbox{sn}(a\omega+b\omega')$, regarded as a function of the elliptic modulus $k^2\; (=x)$.
Therefore, the ansatz for studying the branching of solutions to \eqref{P6} assumed in \cite{guz}
$$
y(x)=\wp\big(\nu_1^{}\omega_1^{}(x)+\nu_2^{}\omega_2^{}(x)+
v(x;\nu_1^{},\nu_2^{});\omega_1^{}(x),\omega_2^{}(x)\big)+\frac{1+x}{3}
\,,
$$
contains a somewhat `redundant' perfect square and can be conveniently replaced with a structure of
Hitchin's type \eqref{hitchin}, \eqref{52}, which simplifies the analysis.

\normalsize
{\em Remark $1$.} In addition to the Picard solution, all its transformations obtained by actions of discrete
symmetries of equation \eqref{P6} \cite[\textsection\,42]{gromak}, for instance, by the variable change
$$
y\mapsto\widetilde y=\frac{y-x}{y-1}=
\frac{\vartheta_4^2(\tau)}{\vartheta_3^2(\tau)}\,
\frac{\theta_3^2\big(A\tau+B|\tau
\big)}{\theta_4^2(A\tau+B|\tau)}\,,
$$
are perfect squares. Adding the modular substitutions from the group
$\boldsymbol{\Gamma}(1)$ and the anharmonic group $S_3$,
transforming the variable $x$ \cite[\textsection\,23]{a}, we obtain the complete set of transformations in \cite[pp.\;214--5]{gromak} and changes of
the parameters
$(\alpha,\beta,\gamma,\delta)$. Together with the Okamoto transformations, all
the transformations obtained above and below expressed in terms of the theta functions and $\theta, \Btau$-functions
reduce to permutations of indices of the functions
$\vartheta,\,\theta,\,\theta'$ in formulae \eqref{simp},
\eqref{hitchin}, \eqref{52}, and \eqref{53}. By
the Picard--Hitchin solution, we therefore understand the whole class of these solutions, without especially
mentioning this.
The Picard--Hitchin solutions are functions of $\sqrt{x}$. On the other hand, the quantity $\sqrt{x}=s$ itself is
related to the second-order linear differential Fuchs equation via the integrals
$K,\,K'(\sqrt{x})$ and
substitution \eqref{subs}. This equation has four singularities located at the points $s=\{0,1,-1,\infty \}$ because the
integrals
$K$ and $K'(x)$ satisfy the Legendre equation with the three singularities
$x_j^{}=\{0,1,\infty \}$. We thus obtain the particular case of the Heun equation
\begin{equation}\label{heun}
Y_{ss}=-\frac12\,\frac{(s^2+1)^2}{s^2\,(s^2-1)^2}\,Y\quad\Rightarrow\quad
\frac{Y_2}{Y_1}=\ri\,\frac{K(s)}{K'(s)}=\tau\,.
\end{equation}
The character of quadrature  \eqref{hitchin}  is therefore described by the following statement.

{\bf Corollary 2.} {\em Painlev\'e equation
\eqref{P6} in the Picard and Hitchin cases is integrable in
$(\theta',\theta)$-function
quadratures over the differential field determined by the Heun functions of form \eqref{heun} or, equivalently, by
the hypergeometric Legendre functions.}

\subsection{One-parameter solutions}
We note that neither formula \eqref{hitchin1} nor the first formula in \eqref{pic} imply
in any obvious way that the second-order poles from the Picard series become the first-order poles for \eqref{52}. Analogously comparing formulae
\eqref{subs}, \eqref{hitchin1}, and \eqref{pic}, we find that the first formula in \eqref{pic} does not provide the
solution behavior in the limit cases for the parameters
$(A,B)$, while this behavior follows automatically
from formula \eqref{52}. For instance, we set  $A\tau+B\to 0$ with $B/A=\alpha$  in \eqref{52} and \eqref{53}. Taking \eqref{weier} and \eqref{EK} into account, we then immediately obtain the $\alpha$-parameter family of solutions
\begin{equation*}
\begin{aligned}
y=\frac{\alpha E'-E+K}{\alpha K'+K}\Leftrightarrow z&=
\wp^{\mbox{-}1}\!\!\left(\frac{\pi\ri}{2}\,\frac{d}{d\tau}\mbox{Ln}
\big\{(\tau+\alpha)\ded^2(\tau)\big\}
\Big|\tau\right)\\[0.3em]
&=\wp^{\mbox{-}1}\!\!\left(\frac{1}{2}\,\frac{\pi\,\ri}{\tau+\alpha}-\eta(\tau)
\Big|\tau\right),
\end{aligned}
\end{equation*}
which can be usefully compared with the corresponding Hitchin formula in \cite[p.\;78]{hitchin}.

Another case looks more complicated, and Hitchin considered it only in the expansion form because
it results in a singular conformal structure. This case is more conveniently analyzed in form \eqref{wpz}, where it
corresponds to the limits
$A\tau+B\to\big\{1,\tau,\tau+1 \big\}$. Setting
$A\tau+B=\varepsilon\,\tau+(1+\varepsilon\alpha)$ and letting
$\varepsilon$  tend to zero, we obtain the solution for which we preserve  the \Btau-form \eqref{entire}:
\begin{equation}\label{L}
\begin{aligned}
\wp(z|\tau)&=\frac{\pi}{2\ri}\frac{d}{d\tau}
\mbox{Ln}\frac{2(\tau+\alpha)\,\dot\vartheta_2+\vartheta_2}
{\ded^2\,\vartheta_2}\\[0.3em]
&=\frac{\pi^2(\vartheta_3^4+\vartheta_4^4)
(\eta'+\alpha\,\eta)+(\tau+\alpha)\big(g_2^{}-\frac12\pi^4
\vartheta_3^4\,\vartheta_4^4
\big)}{\pi^2(\vartheta_3^4+\vartheta_4^4)(\tau+\alpha)+
12\,(\eta'+\alpha\,\eta)}\,.
\end{aligned}
\end{equation}
Here and hereafter, we let the dot above a symbol denote the derivative with respect to $\tau$ and omit the
argument  $\tau$ of the functions $\ded,\,\eta,\,\eta',\,g_2^{},\,
\vartheta$ for brevity. Moreover, we can replace $g_2^{}$ with its
$\vartheta$ equivalent (see appendix). We can calculate this solution in the variables $(x,y)$  if we supplement
relations \eqref{LnK}, \eqref{weier}, \eqref{EK}
with additional relations to obtain the complete basis of the transition between
the
`$\tau$-' and `$x$-representations'. We can do this using the expressions for the quantities $\vartheta$, $\dot\vartheta_2$ and $\eta$:
$$
\vartheta_3^2=\frac{2}{\pi}\,K'\,,\qquad
\dot\vartheta_2=\frac{\ri}{\pi}\,\vartheta_2\,K'E'\,,\qquad
\eta=K'E'-\frac{x+1}{3}\,K'^2 \,.
$$
Simplifying and replacing $\alpha\mapsto \ri\alpha$, we obtain the compact solution
$$
y=\frac{x\,(\alpha K'+K)}{\alpha E'-E+K}\,.
$$
We can analogously simplify the two other limit cases:
$$
y=\frac{2(\alpha K'+K)(E'-K')-\pi}{2(\alpha K'+K)(E'-xK')-\pi}\,x
\quad\Leftrightarrow\quad
\wp(z|\tau)=\frac{\pi}{2\ri}\frac{d}{d\tau}
\mbox{Ln}\frac{2(\tau+\ri\alpha)\,\dot{\vartheta}_3+\vartheta_3}
{\ded^2\,\vartheta_3}\,,
$$
$$
y=\frac{2(\alpha K'+K)(E'-xK')-\pi}{2(\alpha
K'+K)(E'-K')-\pi}\phantom{\,x} \quad\Leftrightarrow\quad
\wp(z|\tau)=\frac{\pi}{2\ri}\frac{d}{d\tau}
\mbox{Ln}\frac{2(\tau+\ri\alpha)\,\dot{\vartheta}_4+\vartheta_4}
{\ded^2\,\vartheta_4}\,.
$$
We see that these solutions are mutually related by transformations of the form
$y\mapsto \frac xy$; they therefore
cannot originate from the modular transformations   $\tau\mapsto \tau+1$ and $\tau\mapsto -\frac{1}{\tau}$.
But we can easily see
that they reflect discrete symmetries of equation
\eqref{th} of the form $z\mapsto \{z+1,z+\tau\}$, and this entails the above
transformation for $y$, which is easily recognizable when written in the standard $\wp(z|\omega,\omega')$-representation:
$$
\wp(z+\omega)-e=\frac{(e-e')(e-e'')}{\wp(z)-e}\quad\Leftrightarrow\quad
\widetilde y= \frac xy\,.
$$
Of course, we can apply this
symmetry to any Picard--Hitchin solution.
We also note that we need the pair of functions $(E,E')$ not
only as a formal completion of the set of functions $(K,K')$:
both the Hitchin and the Chazy solutions
obtained in \cite{mazz} and also any
$\eta,\vartheta,\theta$-function solutions to equation \eqref{P6} are expressed in terms of $(E,E')$.

\subsection{The \protect\Btau-functions}
We also presented examples of the Painlev\'e
$\partial_x\mbox{Ln}$-forms for the degenerate
Hitchin solutions because the \Btau-functions are important for the theory of Painlev\'e equations, being the
generators of their Hamiltonians and other objects
\cite{iwasaki,conte,okamoto,okamoto2}. But the  \Btau-functions are primarily
important not for constructing Hamiltonians, which are known for all Painlev\'e equations, but for representing
solutions in form \eqref{entire} see \cite{okamoto2}
and the explanations in \cite[pp.\;740--1]{conte} in regard to constructing
\Btau-functions
corresponding to the positive and negative residue series for solutions of the equation $\mathcal{P}_{2}$.

We note that other known examples (see the citations at the end of sect.~2.1) require the existence of the explicit transition ${}_2F_1\mapsto \tau$ because the ${}_2F_1$-series are essentially nonholomorphic functions with branching.
Reducing them to the
`genuine' $(\tau)$-form is a very complicated problem, which is equivalent to inverting
the integral ratios for equations of the Fuchsian type. In this respect, formulae are always `final' in the
uniformizing $\tau$-representation
(sect.~5), while the transition ${}_2F_1\mapsto \tau$ is nontrivial and may contain
several hypergeometric series. For instance, in the above examples, besides the integrals $K,K',E,E'$, \ie\ functions of the type
${}_2F_1\big({\pm}\frac12,\frac12;1\big|z\big)$, we encounter the function $\ded^2(\tau)$. It can be shown \cite[p.\;263]{br} that the
hypergeometric equation related to this function is just the equation for the function
${}_2F_1\big(\frac{1}{12},\frac{1}{12};\frac23\big|J\big)$, \ie\
the equations that determine the classic Klein modular invariant
$J=(g_2^3/\Delta)(\tau)$ \cite[\S\,10]{a},
\cite[\S\,14.6.2]{bateman}:
$$
J(J-1)\frac{d^2}{dJ^2}\ded^2+\frac16(7J-4)\frac{d}{dJ}\ded^2+
\frac{1}{144}\ded^2=0\,.
$$

By this means we currently have only the Hitchin class in which all the results including degenerations
can be written explicitly. It remains to note that the found Painlev\'e form admits a natural interpretation,
analogous to formulae of the straight-line sections of the two-dimensional Jacobians in soliton equation
theory. Function \eqref{mero} is the canonically normalized meromorphic elliptic integral
$I(z-B|\tau)$, which has
only the first-order pole at the point $z=B$
and is considered on the straight-line section $z-A\tau+B=0$
of the space $\{\mathbb{C}\times\mathbb{H}^+\}$, obtained as the direct product of the elliptic curve (the parameter $z$) and the
one-dimensional moduli space of elliptic curves (the parameter $\tau$).

The \Btau-function constructed in \cite{korotkin}, corresponds only to Picard poles \eqref{poles}  and coincides with the function
$$
\Btau_{\!\!1}^{}(x;A,B)=\theta_1 \!\!\left(\!
A\mfrac{K(\sqrt{x})}{K'(\sqrt{x})}+B\Big|\mfrac{\ri\,
K(\sqrt{x})}{K'(\sqrt{x})}\right) \,.
$$
up to a fixed nonmeromorphic factor. Hence, there exists a second Hamiltonian corresponding to the
second
\Btau-function and poles with the opposite signs of residues. This function is given by the formula
$$
\Btau_{\!\!2}^{}(x;A,B)= \Btau_{\!\!1}^{}(x;A,B)\,\frac{d}{dB}
\mbox{Ln}\,\big\{ \Btau_{\!\!1}^{}(x;A,B)\,\re^{2\pi A
B}_{\mathstrut}\big\}\,.
$$
We choose such a representation because meromorphic Abelian integrals (integrals of the second kind) of
the type of functions \eqref{mero}
can be represented as derivatives of logarithmic integrals (integrals of the third
kind) with respect to the parameter determining the location of one of the two logarithmic poles (with
respect to the constant $B$).
This may entail further generalizations, but we do not discuss this here.

\section{The equation $\mathcal{P}_{6}$ and uniformization of curves}

\subsection{Algebraic solutions} Having the general integral leads to several corollaries, for instance, to
the limit cases with respect to the parameters $(A,B)$, for example,

\begin{equation}\label{AB}
A\,\tau+B=\frac nN\,\tau+\frac mN
\end{equation}
We thus obtain an infinite series of algebraic solutions of the form $P(x,y)=0$, where $P$  is a polynomial in
its arguments. Hitchin also presented particular solutions of the algebraic type, and the origin of the whole
algebraic Picard series was revealed in \cite{mazz} where asymptotic regimes and monodromies related to this and
other solutions and the character of their branching were also considered.

It is less known that some of these results and also the very scheme for constructing such solutions were
already contained in Fuchs' paper [24] with the same title as his well-known paper\footnote{In a one-line footnote in \cite[p.\;300]{picard} Picard also claimed that the infinite series of algebraic solutions appears in the case where the
numbers $(A,B)$ are commensurable.} of 1907. The procedure for constructing solutions was missed in \cite{fuchs2} but the method for obtaining these solutions is obviously
the same: the theorems of multiplication and addition for elliptic functions can be used. Moreover, Fuchs
pointed out \cite{fuchs2} how to use the known Kiepert determinant multiplication formula for $\wp(Nz)$ as a rational
function of $\wp(z)$ \cite[p.\;19]{we2},
\cite[p.\;332]{WW}; he also wrote the Puiseux series for these solutions \cite[\textsection\,IV]{fuchs2} and considered
the monodromy of the Fuchs equations in great detail. Such determinant schemes can be interesting in
themselves because the list of algebraic solutions of the Painlev\'e equations steadily increases. We restrict
ourself to the statement (without proof) that we can obtain effective recursion schemes for constructing
such solutions. We do not present examples of these recursions because another immediate corollary from
what was said above seems more interesting to us.

From the explanations in sect.~1.1 and from explicit form
\eqref{53} of the $z(\tau)$-dependence for the solution to equation \eqref{P6wp}, we obtain infinite series of algebraic curves $P(u,v)=0$
of nontrivial genera $g>1$ with the explicit
parameterizations by the single-valued analytic functions $u(\tau),\,v(\tau)$. Because the Kiepert formula above
contains a redundant set of functions
$(\wp,\wpp,{\textstyle\wp\mbox{\small$''$}},\ldots)$, it is more convenient to pass to the language of theta
functions, which also provides a more compact answers. In this case, it is desirable to have multiplication
theorems directly for the functions $\theta_k$, $\dtheta$.
To the best of our knowledge, such formulae are absent in
the literature, and we therefore present them in the appendix (Lemma\;2). For the Hitchin solutions, we
can use transformation \eqref{ok}, but we lack the efficiency when acting with them on algebraic functions whose
complexity rapidly increases as $N$ increases.
It is therefore more convenient to use the abovementioned
theorems for the
$\theta,\theta'$-functions directly in formulae \eqref{Tau}, \eqref{53} (also see Proposition \;2 below). It is
obvious how we can apply these theorems to the solutions under investigation, and we therefore pass to demonstrating specific examples and their corollaries.

\subsection{Examples}
Setting $A\tau+B=\frac13$, we reconstruct the Dubrovin solution
\cite[(4.5)]{mazz}
\begin{equation}\label{pic2}
y^4-(6\,y-4)\,y\,x +(4\,y-3)\,x^2=0\,.
\end{equation}
Formula \eqref{subs} then results in the curve of genus $g=3$ of the form
\begin{equation}\label{curve1}
4\,(u^2+u^{-2})\,v=v^4-6\,v^2-3\,,\qquad
\left\{u=\frac{\vartheta_4(\tau)}{\vartheta_3(\tau)}\,,\quad v=
\frac{\theta_2^2\big(\frac16\big|\tau\big)}
{\theta_1^2\big(\frac16\big|\tau\big)}\right\}
\end{equation}
if we use the substitution $(x=u^4,\; y =-u^2\,v)$. This curve does not belong to the series of classic modular
equations
\cite[\textsection\,14.6.4]{bateman} nor to their closest variations
\cite[\textsection\,7]{br}. In turn, it is remarkable that because of the birational transformation
$$
\left\{\!\!\begin{array}{rcl} \ds
u\!\!\!&=&\ds\!\!\!\frac{\boldsymbol{x}^4-\boldsymbol{y}-1}{4\,\boldsymbol{x}}_{{}_{\ds\mathstrut}}\\
\ds
v\!\!\!&=&\ds\!\!\!\frac{\boldsymbol{x}^4-\boldsymbol{y}+1}{2\,\boldsymbol{x}^2}
\end{array},\qquad\quad\right\{\!\!
\begin{array}{rcl}
\ds \boldsymbol{x}\!\!\!&=&\ds\!\!\!\frac{u\,(v^2-1)}{2\,(u^2+v)}_{{}_{\ds\mathstrut}}\\
\ds \boldsymbol{y} \!\!\!&=&\ds\!\!\!
\frac{(v^2+3)^2}{v^2\,(v^2-1)}\!\left(\frac{4\,u^2\,v}{v^2+3}+1\right)
\end{array},
$$
equation \eqref{curve1} is isomorphic to the classic hyperelliptic Schwarz curve
\begin{equation}\label{x8}
\boldsymbol{y}^2=\boldsymbol{x}^8+14\,\boldsymbol{x}^4+1\,.
\end{equation}
In Schwarz's Gesammelte Abhandlungen \cite{schwarz} this curve appears in many contexts in both volumes\footnote{For instance, it appears in the 1867 paper  \cite[{\bf I}: p.\;13]{schwarz}, where the famous Schwarz derivative $\boldsymbol{\mathit{\Psi}}(s,u)$ was introduced.}, but no variants of its
parametrization were known. We present the formula for $\boldsymbol{x}(\tau)$,
because it admits the compact simplification
\begin{equation}\label{xuni}
\boldsymbol{x}=\frac{\ded^3(\tau)\,\theta_2\big(\frac13\big|\tau\big)}
{\theta_1^2\big(\frac16\big|\tau\big)\,
\theta_3^2\big(\frac16\big|\tau\big)}\,.
\end{equation}

For this curve and other curves and functions introduced by this method, we can write the remaining
uniformization attributes (which is a purely technical problem): differential equations for function $y(\tau)$, the
${}_2F_1$-solutions of the corresponding Fuchsian equations, the representations of the monodromy group of
these equations, the Poincar\'e polygons
\cite{bateman}, the Abelian integrals as functions of $\tau$, their reductions, if any, to elliptic integrals, etc.

If we consider analogous solution \eqref{hitchin1}, with the same parameters $A\tau+B=\frac13$, then we obtain the
curve that is the Hitchin solution
\begin{equation}\label{hh}
y^4-2\,(x+1)\,(y^2+x)\,y+6\,x\,y^2+x^3-x^2+x=0\,.
\end{equation}
But this curve with the substitution $x=u^4$ must become isomorphic to the previously described curves
because the initial Picard \eqref{pic2} and Hitchin \eqref{hh} solutions have the genus $g=0$. In turn, formulae
\eqref{subs} and \eqref{hitchin1} provide an equivalent but not obvious parameterization of \eqref{hh}
$$
y=\frac{\vartheta_4^2(\tau)} {\vartheta_3^2(\tau)}\! \left(
\!\pi\,\vartheta_2^2(\tau)\,\frac{\theta_3\,\theta_4}{\theta_2\,\dtheta}
\big(\tfrac16\big|\tau\big) -1\right)
\frac{\theta_2^2}{\theta_1^2}\big(\tfrac16\big|\tau\big)\,,
$$
which follows from the nontrivial
$\vartheta,\,\theta,\,\theta'$-function relations\footnote{In particular, $\dtheta(p\,\tau\!+\!q|\tau)$ is algebraically related to the functions $\theta(p\,\tau\!+\!q|\tau)$, for rational $p,q$ while this relation can
be only differential in the generic case (see Lemma~1 in the appendix).}. The general statement about the algebraic
Hitchin solutions follows from formulae
\eqref{subs} and \eqref{simp}.

{\bf Proposition 2.} {\em The algebraic Hitchin solutions admit the parameterization}
\begin{equation}\label{param}
x=\frac{\vartheta_4^4}{\vartheta_3^4}\,,     \qquad
y=\frac{\vartheta_4^2}{\vartheta_3^2}\,\frac{\theta_2}{\theta_1^2}
\left\{
\frac{\pi\,\vartheta_2^2\,\theta_3\,\theta_4}{\dtheta+2\pi\ri\frac
nN\theta_1}-\theta_2 \right\},
\end{equation}
{\em  where the $\theta,\theta'$-functions are taken with the arguments in \eqref{AB}.
Functions  \eqref{param} satisfy autonomous third order
differential equations of the forms
$$
\frac{\dddot {x\mspace{2mu}}}{\dot x^3}-\frac32\,\frac{\ddot
x^2}{\dot x^4}=-\frac12\,\frac{x^2-x+1}{x^2(x-1)^2}\,,\qquad
\frac{\dddot y}{\dot y^3}-\frac32\,\frac{\ddot y^2}{\dot
y^4}=\boldsymbol{\mathcal{Q}}(x,y)
$$
with the explicitly calculable rational function
$\boldsymbol{\mathcal{Q}}(x,y)$ for all $(n,m,N)$.}

Because of this proposition and the simplicity of formulae \eqref{param}, there is no need to use the equations  $P(x,y)=0$ because the uniformization allows investigating solutions with arbitrarily large
 $(n,m,N)$. Here, we
include differentiation, the Laurent--Puiseux series, calculation of monodromies, the branching, the genus,
plotting of graphs, etc.

The zero genus is not a general property of the solutions under study. For example, under
$A\tau+B=\frac15$  the Picard solution has the form of a cumbersome but elliptic curve
$P(\,\mbox{\tiny\raisebox{0.85em}{6}\hspace{-0.65em}}{x},
\mbox{\tiny\raisebox{0.85em}{12}\hspace{-0.85em}}{y}\,)=0$:
\begin{equation}\label{15}
\begin{aligned}
y^{12}&-50xy^{10}+140(x^2+x)y^9-5(32x^3+89x^2+32x)y^8\\
&+16(4x^4+35x^3+35x^2+4x)y^7-
60(4x^4+13x^3+4x^2)y^6+360(x^4+x^3)y^5\\
&-105x^4y^4-80(x^5+x^4)y^3+
2(8x^6+47x^5+8x^4)y^2-20(x^6+x^5)y+5x^6=0
\end{aligned}
\end{equation}
with the Klein modular invariant $J=\frac{256}{135}$. The Weierstrass $\wp$-parameterization for \eqref{15} is very complicated and less adjustable for this curve, while for any $\tau\in\mathbb{H}^+$ we have
$$
x=\frac{\vartheta_4^4(\tau)} {\vartheta_3^4(\tau)}\,,\qquad
y=-\frac{\vartheta_4^2(\tau)}{\vartheta_3^2(\tau)}\,
\frac{\theta_2^2\big(\frac{1}{10}|\tau
\big)}{\theta_1^2\big(\frac{1}{10}|\tau \big)}\,.
$$
The curve genus increases rapidly.  For $N=4,5,6,7,8,9,\ldots$ the genera of the respective Picard--Hitchin
solutions are $g=0,1,1,4,5,7,\ldots$. We mention here that the genus depends on the choice of
the representing equation with the Painlev\'e property. For instance, the same solution \eqref{15}, but in Picard
representation (\ref{su}--\ref{u}) has a nonhyperelliptic genus $g=3$, and a simple solution of equation \eqref{u}, which is
birationally isomorphic to the lemniscate $w^2=z^4-1$ (the calculation is nontrivial),
already appears at $A\tau+B=\frac14$:
\begin{equation}\label{xu}
\begin{aligned}
x^4\,u^{16}-20\,x^3\,u^{12}+32\,x^2(x+1)\,u^{10}
-2\,x\,(8\,x^2+29\,x+8)\,u^8&\\
{}+32\,x\,(x+1)\,u^6-20\,x\,u^4+1&=0\,.
\end{aligned}
\end{equation}

{\em Remark $2$.} In \cite{br}, we advanced an assumption about a possible relation between the theta constants $\theta\big(u(\tau)|\tau\big)$ of general form and the uniformizing Fuchs equations. A particular case of such a constant is just
the Dedekind function
$$
\ded(\tau)=-\ri\,\re^{\frac{\pi\ri}{3}\tau}_{\mathstrut}\,
\theta_1(\tau|3\tau)\,,
$$
which appears both in our formulae and in some known but
isolated parameterizations of modular equations (mostly of the genus zero; see the bibliography in \cite{br}). It
follows from the above results that the Painlev\'e equation admits an infinite class of generalizations in this
direction using the fundamental object $\theta\big(\frac
nN\tau+\frac mN|\tau\big)$, together with the explicit formulae for curves that
are not necessarily solutions of the equation $\mathcal{P}_{6}$, for functions on these curves, for the differential equations
on these functions corresponding to the subgroups of $\boldsymbol{\Gamma}(2)$, and so on.

We note the following computational fact. The Painlev\'e curves are specific and not always easily
accessible by various algorithms. Curve \eqref{xu} is a good example. But in the framework of what was said
above, convenient equations appear, which in particular include an important class of hyperelliptic curves.
Even less obvious is that function \eqref{xuni} in fact corresponds to a larger curve
\begin{equation}\label{x85}
\boldsymbol{z}^2= (\boldsymbol{x}^8+14\,\boldsymbol{x}^4+1)
(\boldsymbol{x}^5-\boldsymbol{x})\,.
\end{equation}
Namely, function \eqref{xuni}  also parameterizes this equation in the sense that $\boldsymbol{z}(\tau)$, determined from Eqs. \eqref{xuni} and
\eqref{x85} is an analytic function that is uniquely defined in its whole domain, \ie, in $\mathbb{H}^+$. This example
provides one more realization of a rather nontrivial tower of embedded hyperelliptic curves with the universal
uniformizing function \mbox{\protect\raisebox{0.2em}{$\chi$}}$(\tau)$ (see \cite{br} about this function).

We omit both the derivation of formula \eqref{x85} and other theta-functional consequences of what was said
above. Because of their large number, they deserve a separate investigation that is beyond the scope of this
paper. This material is presented in work \cite{br3}.

\section{Concluding remarks}
\noindent Using the analytic results in both
$x$-representation \eqref{hitchin}, \eqref{52}, and
$\tau$-representation \eqref{Tau}, \eqref{53} we can
easily obtain several useful formulae, for example, simple expressions for the conformal factor $F$ of the
Tod--Hitchin metric
\cite{babich,hitchin,tod}. This factor was expressed in terms of the derivatives of the theta functions
with respect to the integration constants in \cite[theorem~5.1]{babich},  and the solution $y(x)$ was expressed in terms of the double
derivatives of the theta functions with respect to both its arguments. By Lemma\;1 the final
answers are actually expressed in the basis
$\vartheta,\dtheta,\theta_k(A\tau+B|\tau)$. Another example is third-order differential
equations for the \Btau-functions and for the entire holomorphic functions whose ratio gives the solution $y(x)$. On numerous occasions, Painlev\'e himself indicated the existence and importance of these equations (see,
e.g., \cite[p.\;1114]{painleve}). They follow from the third-order equations for $\sigma(a\tau+b|\tau)$,
which follow from Lemma\;3 in the appendix.

Recalling the remark in sect~1.1, we note that the known solutions and local branchings
\cite{gromak,mazz,hitchin,guz} are compatible with the uniformization by the Painlev\'e substitution, \ie, they generate the
$\tau$-representations
of the solutions in terms of the single-valued functions $y(\tau)$. The necessary data for establishing this property
for all the solutions is presumably contained in Theorem\;2 in \cite{okamoto2} and in results of work by Guzzetti \cite{guz}. But in the uniformization
context, it is desirable to have the direct statement about the global single-valuedness of the function $\wp\big(z(\tau)|\tau\big)$ for all values of $(\alpha,\beta,\gamma,\delta)$. This would provide a new important interpretation of the equation $\mathcal{P}_{6}$. Equation \eqref{P6} would differ from the lower Painlev\'e equations in that their solutions are globally meromorphic
on the plane $\mathbb{C}$; for the equation
$\mathcal{P}_{6}$, they are globally meromorphic on the upper half-plane $\mathbb{H}^{+}$ with
the constant negative Lobachevskii curvature, which is the universal covering of $\overline{\mathbb{C}}\backslash \{0,1,\infty\}$. A broad
class of algebraic solutions is then provided by automorphic functions on the Riemann surfaces with three
punctures (orbifolds), while the transcendental solutions of the Picard--Hitchin type themselves generate
another, physically meaningful class of single-valued functions. In this context, solutions to equation \eqref{P6wp} do not
follow automatically from those of
\eqref{P6}, and integrating
\eqref{P6wp} is therefore an independent problem.

\section{Appendix: rules for differentiating the theta-function}
\subsection{Jacobi's $\theta$-functions} These functions are used in the following definitions
\cite{a}:
\begin{alignat*}{6}
\theta_1(z|\tau)&= -\ri\,\re^{\frac14\pi\ri\tau}_{\mathstrut}
\sum\limits_{\sss k=-\infty}^{\sss \infty}\! (-1)^k\,
\re^{(k^2+k)\pi\ri\,\tau}\, \re^{(2k+1)\pi \ri\, z}\,,&\qquad&
\theta_3(z|\tau)= \sum\limits_{\sss k=-\infty}^{\sss \infty}
\re^{k^2\pi\ri\,\tau}\,\re^{2k\pi \ri\,
z}\,,\\
\theta_2(z|\tau)&=
\phantom{-\ri}\,\re^{\frac14\pi\ri\tau}_{\mathstrut}
\sum\limits_{\sss k=-\infty}^{\sss \infty}
\re^{(k^2+k)\pi\ri\,\tau}\,\re^{(2k+1)\pi \ri\,
z}\,,
&&
\theta_4(z|\tau)= \sum\limits_{\sss
k=-\infty}^{\sss \infty}\! (-1)^k\, \re^{k^2\pi\ri\,\tau}\,
\re^{2k\pi \ri\, z}\,.
\end{alignat*}
We also introduce the fifth independent object
$$
\theta'_1(z|\tau)= \pi\,\re^{\frac14\pi\ri\tau}_{\mathstrut}
\sum\limits_{\sss k=-\infty}^{\sss \infty}\! (-1)^k\,
(2k+1)\,\re^{(k^2+k)\pi\ri\,\tau}\, \re^{(2k+1)\pi \ri\, z}
$$
and the $\vartheta$-constants are defined as $\vartheta_k \deq
\vartheta_k(\tau)=\theta_k(0|\tau)$ and $\vartheta_1\equiv 0$.

{\bf Lemma 1.} {\em The five Jacobi functions
$\theta_1\,,\,\theta_2\,,\,\theta_3\,,\,\theta_4$ and $\,\dtheta\deq
\mfrac{\partial\theta_1}{\partial z}$ satisfy the closed ordinary differential
equations with respect to the variables $(z,\tau)$  over the field of coefficients $\eta(\tau)$ and $\vartheta^2(\tau)$}:
\begin{align}
& \label{DthetaX}
\left\{
\begin{aligned}
\frac{\partial\theta_k}{\partial z}&= \frac{\dtheta}{\theta_1}
\,\theta_k- \pi\,\vartheta_k^2\!\cdot\!
\frac{\theta_\nu\,\theta_\mu}{\theta_1} \\[0.3em]
\ds\mspace{1mu}\frac{\partial\dtheta}{\partial z}&=
\frac{\dtheta^2}{\theta_1}-\pi^2\vartheta_3^2\,\vartheta_4^2
\!\cdot\! \frac{\theta_2^2}{\theta_1}-
4\!\left\{\eta+\frac{\pi^2}{12}\big(\vartheta_3^4+\vartheta_4^4\big)
\right\} \!\cdot\!\theta_1
\end{aligned}\right.,
\\[0.3em]\notag
& \label{tau} \left\{
\begin{aligned}
\frac{\partial\theta_k}{\partial \tau}&= \frac{-\ri}{4\,\pi}\,
\frac{\dtheta{}^2}{\theta_1^2} \,\theta_k+ \frac{\ri}{2}\,
\vartheta_k^2\!\cdot\!\dtheta\,
\frac{\theta_\nu\,\theta_\mu}{\theta_1^2}+ \frac{\pi\ri}{4}\,\Big\{
\vartheta_3^2\,\vartheta_4^2\!\cdot\!\theta_2^2-
\vartheta_k^2\,\vartheta_\mu^2\!\cdot\!\theta_\nu^2-
\vartheta_k^2\,\vartheta_\nu^2\!\cdot\!\theta_\mu^2
\Big\}\,\frac{\theta_k}{\theta_1^2}\\[0.3em]
&\quad +\frac{\ri}{\pi}\!\left\{ \eta+\frac{\pi^2}{12}
\big(\vartheta_3^4+ \vartheta_4^4\big)\right\}\!\cdot\!
\theta_k\\[0.3em]
\frac{\partial \dtheta}{\partial\tau}&=
\frac{-\ri}{4\,\pi}\,\frac{\dtheta{}^3}{\theta_1^2}
+\frac{3\,\ri}{\pi}\!\left\{
\frac{\pi^2}{4}\vartheta_3^2\,\vartheta_4^2\!\cdot\!
\frac{\theta_2^2}{\theta_1^2} +\eta+
\frac{\pi^2}{12}\big(\vartheta_3^4+\vartheta_4^4\big)
\right\}\dtheta -\frac{\pi^2}{2}\ri\,\vartheta_2^2\,\vartheta_3^2\,
\vartheta_4^2\!\cdot\!
\frac{\theta_2\theta_3\theta_4}{\theta_1^2}
\end{aligned}\right.\,,
\end{align}
{\em where}
$$
\nu=\frac{8\,k-28}{3\,k-10}\,,\qquad
\mu=\frac{10\,k-28}{3\,k-8}\,,\qquad k=1,2,3,4\,.
$$

These differentiations are not completely symmetric, and the equations themselves mean that not only
elliptic functions (the ratios
$\theta_j/\theta_k$) have the closed differential calculus with respect to both the variables $(z,\tau)$ independently but also the same holds for their `constituents', the theta functions including the
function $\dtheta$. The details and, not less important, the analysis of the integrability condition for the above
equations have been detailed in work \cite{br2}. The basis
$\theta_{1,2,3,4}$ can be obviously closed by adding any of the functions $\theta_k'$, not necessarily $\dtheta$. It is essential that the known equation $4\pi\ri\,\theta_\tau=\theta_{\mathit{zz}}$ must be treated as a corollary of the
above equations and not vice versa because the Lemma~1 defines ordinary differential equations,
while this equation is an equation of the heat kernel type in partial derivatives; the theta functions are just particular solutions of it.

{\bf Lemma 2.} {\em Let $n$ be an integer and $n_1^{}\deq n-1$. Then the Jacobi functions satisfy the recursive
multiplication formulae}
\begin{equation}\label{prime1}
\left\{
\begin{array}{l}\ds
\theta_1(2z)=2\,\theta_1(z)\,
\frac{\theta_2(z)\,\theta_3(z)\,\theta_4(z)}
{\vartheta_2\,\vartheta_3\,\vartheta_4}\\\\
\ds \theta_1(nz)=\frac{\theta_3^2(n_1^{} z)\,\theta_2^2(z)-
\theta_2^2(n_1^{}z)\,\theta_3^2(z)} {\vartheta_4^2 \cdot
\theta_1\big((n-2)z\big)}
\end{array}\right.\,,
\end{equation}
$$
\left\{
\begin{array}{l}
\ds \theta_2^{}(nz)=\frac{\theta_3^2(n_1^{}z)\,\theta_3^2(z)-
\theta_4^2(n_1^{}z)\,\theta_4^2(z)} {\vartheta_2^2 \cdot
\theta_2^{}\big((n-2)z\big)}\\
\ds
\theta_3^{}(nz)=\frac{\theta_2^2(n_1^{\ds\mathstrut}z)\,\theta_2^2(z)+
\theta_4^2(n_1^{}z)\,\theta_4^2(z)} {\vartheta_3^2 \cdot
\theta_3^{}\big((n-2)z\big)}\\
\ds
\theta_4^{}(nz)=\frac{\theta_3^2(n_1^{\ds\mathstrut}z)\,\theta_3^2(z)-
\theta_2^2(n_1^{}z)\,\theta_2^2(z)} {\vartheta_4^2 \cdot
\theta_4^{}\big((n-2)z\big)}
\end{array}\right.\,.
$$
{\em We obtain the multiplication formula for the function $\dtheta(nz)$ by taking the derivative in \eqref{prime1}, and subsequently
using formulae \eqref{DthetaX}. The lemma statement also holds for an arbitrary complex $n$.}

It is natural here to assume that there are general nonrecursive expressions for these formulae in form of certain
determinants.

\subsection{The Weierstrass functions}
The modular Weierstrass functions $g_2^{}(\tau)$,
$g_3^{}(\tau)$ and $\eta(\tau)$ are
defined by the series
\begin{alignat}{2}\notag
\eta(\tau)&=\phantom{0}2\,\pi^2 \left\{ \frac{1}{24}-
\sum\limits_{k=1}^\infty\,
\frac{\re^{2k\pi\ri\tau}}{\big(1-\re^{2k\pi\ri\tau}\big)^2}
\right\},\\[0.3em]\label{g23}
g_2^{}(\tau)&= 20\,\pi^4 \left\{\frac{1}{240} +
\sum\limits_{k=1}^\infty\,
\frac{k^3\,\re^{2k\pi\ri\tau}}{1-\re^{2k\pi\ri\tau}}
\right\}_{\ds{}_{\mathstrut}}=
\frac{\pi^4}{24}\big(\vartheta_2^8+\vartheta_3^8+\vartheta_4^8\big)\,,
\\[0.3em]\notag
g_3^{}(\tau) &= \;\frac73\,\pi^6 \left\{\frac{1}{504} -
\sum\limits_{k=1}^\infty\,
\frac{k^5\,\re^{2k\pi\ri\tau}_{}}{1-\re^{2k\pi\ri\tau}} \right\}=
\frac{\pi^6}{432}\big(\vartheta_2^4+\vartheta_3^4\big)
\big(\vartheta_3^4+\vartheta_4^4\big)
\big(\vartheta_4^4-\vartheta_2^4\big)\,.
\end{alignat}

{\bf Lemma 3.} {\em The Weierstrass functions
$\sigma,\zeta,\wp,\wpp(z|\tau)$ satisfy the nonautonomous dynamical
system with the coefficients $\eta(\tau)$, $g_{2}^{}(\tau)$
and with the parameter $z$\/}:
\begin{equation}\label{W}
\left\{
\begin{array}{l}
\ds\;\frac{\partial\sigma}{\partial\tau}= \phantom{-}\frac{\ri}{\pi}
\Big\{ \wp-\zeta^2
+2\,\eta\,(z\,\zeta-1)-\Mfrac{1}{12}\,g_2^{}\,z^2\Big\}\,\sigma\\[0.8em]
\left. \!\!\!
\begin{array}{l}
\ds\;\frac{\partial\zeta}{\partial\tau}= \phantom{-} \frac{\ri}{\pi}
\Big\{ \wpp+2\,(\zeta-z\,\eta)\,\wp+ 2\,\eta\,\zeta
-\Mfrac16\,g_2^{}\,z\Big\}\\[0.8em]
\ds\;\frac{\partial\wp}{\partial\tau}=-\frac{\ri}{\pi} \Big\{
2\,(\zeta-z\,\eta)\,\wpp+4\,(\wp-\eta)\,\wp
-\Mfrac23\,g_2^{}\Big\}\\[0.8em]
\ds\frac{\partial\wpp}{\partial\tau}= -\frac{\ri}{\pi}\mspace{1.5mu}
\bbig\{ 6\,(\wp-\eta)\,\wpp+(\zeta-z\,\eta)(12\,\wp^2-g_2^{})
\bbig\}
\end{array}\right\}
\end{array}
\right.,
\end{equation}
{\em where we use the right brace additionally to denote the differential closedness of the functions $\zeta,\wp,\wpp$. Using the replacements
$\zeta\mapsto\zeta-z\eta$ and
$\sigma\mapsto\sigma\exp\!\big\{{-}\frac12\eta z^2\big\}$ we eliminate the parameter $z$ from Eqs. \eqref{W}, which is equivalent to passing to the functions $\theta'$ and $\theta$.}

Because $\wpp$ can be expressed algebraically in terms of $\wp$, the `essential' part of system \eqref{W} is its
second and third equations. In turn, they are equivalent to the equations for the functions
$\zeta(A\tau+B|\tau)$ and $\wp(A\tau+B|\tau)$, if we take the known relation $\partial_z\zeta(z|\tau)=-\wp(z|\tau)$ into account. Therefore, if we regard
the parameters $(\alpha,\beta,\gamma,\delta)$ as the moduli of the equation $\mathcal{P}_{6}$, then the Picard--Hitchin--Okamoto class
(see Remark~1) has  common moduli, and all the equations $\mathcal{P}_{6}$ in this class are then equivalent to a single
representative, which is system \eqref{W}.

{\bf Corollary 3.} {\em The canonical representative of the equations
$\mathcal{P}_{6}$ in the class of the Picard--Hitchin--Okamoto parameters is the system of two equations for the functions \mbox{$\Z$, $\WP$}}:
\begin{equation}\label{canon}
\frac{d\Z}{d\tau}=\frac{\ri}{\pi} \Big\{ \WPP+
2\,(\WP+\eta)\,\Z\Big\}\,,\qquad \frac{d\WP}{d\tau}=-\frac{\ri}{\pi}
\Big\{ 2\,\Z\,\WPP+4\,(\WP-\eta)\,\WP -\mfrac23\,g_2^{}\Big\}
\end{equation}
{\em and its common integral}
$$
\Z=\zeta(A\tau+B|\tau)-A\eta'(\tau)-B\eta(\tau)\,,\qquad
\WP=\wp(A\tau+B|\tau)\,.
$$
{\em Here \WPP\  $=\sqrt{4\WP^3-g_2^{}\WP-g_3^{}}$, and Eqs. \eqref{canon} are supplemented by their corollary\/}:
$$
\frac{d\WPP}{d\tau}= -\frac{\ri}{\pi}\, \bbig\{
6\,(\WP-\eta)\,\WPP+(12\,\WP^2-g_2^{})\,\Z \bbig\}\,.
$$

It is clear that we can keep either of the quantities $\WP,\,\WPP$, and solutions and their derivatives are
then rational functions of $(\Z,\WP,\WPP)$.
Such derivatives can then be treated as the $\wp$-form of Okamoto
transformations \eqref{ok}. For instance, for the Hitchin solution
$H=\wp(z|\tau)$, we have
$$H=\WP+\frac12\,\frac{\WPP}{\Z}\,,\qquad \frac{\pi}{\ri}\,\frac{d
H}{d\tau}=
\frac{4\,H^3-g_2^{}\,H-g_3^{}}{\WP-H}+2\,H^2+4\,\eta\,H+\frac16\,g_2^{}\,.
$$
Choosing a convenient basis composed of $(\Z,\WP,\WPP)$ based on governing equations \eqref{canon}, we can
go further and develop the inverse, \ie, the integral Painlev\'e calculus in both the uniformizing
$(\tau)$-half-plane
and the $(x)$-plane. There are numerous examples, and in addition to the first obvious example following from
the
\Btau-form of \eqref{hitchin}
$$
\int\limits^{\;x}\!\!\frac{y\,dx}{x(x-1)}= -2\,\mathrm{Ln}
\mfrac{\mbox{\normalsize$\dtheta$}\!\bbig(
A\frac{K}{K'}+B\bbig|\frac{\ri\, K}{K'}\bbig)+
\mbox{\normalsize$2\,\pi\, A\!\cdot\!\theta_1$}\!\bbig(
A\frac{K}{K'}+B\bbig|\frac{\ri\, K}{K'}\bbig) }
{\mbox{\normalsize$\sqrt{1-x}\,\,K'\cdot\theta_1$} \!\bbig(
A\frac{K}{K'}+B\bbig|\frac{\ri\, K}{K'}\bbig)}\,,
$$
we present another, more remarkable example:
$$
\frac{\ri}{\pi}\int\limits^{\;\tau}\!\!\big(\WP-\Z^2
\big)d\tau=\mbox{Ln}
\,\theta_1\!\big({\textstyle\frac12}A\tau+{\textstyle\frac12}B|\tau\big)-
\mbox{Ln}\,\ded(\tau)+\frac{\pi\ri}{4}\,A^2\tau\,.
$$
This example, which generates all other integrals, can be logically considered the Painlev\'e $(\tau)$-analogue of the
well-known integral Weierstrass relation between the meromorphic objects $\zeta,\,\wp$ and the entire function $\sigma$. Again imposing condition
\eqref{AB}, we obtain various theta constants $\theta\big(u(\tau)|\tau\big)$,
$\theta'\big(u(\tau)|\tau\big)$ of the general
form in the left-hand side and the value of the integral calculated on these constants, itself a theta constant,
in the right-hand side.

\end{document}